\documentclass[11pt]{amsart}

\usepackage{amssymb,amsmath}
\usepackage{amsthm}
\usepackage{qsymbols}

\newtheorem{thm}[equation]{Theorem}
\newtheorem{lem}[equation]{Lemma}
\newtheorem{prop}[equation]{Proposition}
\newtheorem{cor}[equation]{Corollary}
\newtheorem{defi}[equation]{Definition}
\newtheorem{rem}[equation]{Remark}

\newcommand{\sd}{\prec_{\mathrm{st}}}

\newcommand{\st}{\text{ such that }}
\newcommand{\rf}{\mathrm{ref}}

\newcommand{\la}{\left\{}
\newcommand{\ra}{\right\}}
\newcommand{\lp}{\left(}
\newcommand{\rp}{\right)}
\newcommand{\lc}{\left[}
\newcommand{\rc}{\right]}

\newcommand{\ep}{\varepsilon}

\newcommand{\R}{\mathbb{R}}
\newcommand{\N}{\mathbb{N}}

\def\AArm{\fam0 \rm}%
\newdimen\AAdi%
\newbox\AAbo%
\def\AAk#1#2{\setbox\AAbo=\hbox{#2}\AAdi=\wd\AAbo\kern#1\AAdi{}}%
\newcommand{\1}{{\ensuremath{{\AArm 1\AAk{-.8}{I}I}}}}

\newcommand{\E}{\mathbb{E}}

\newcommand{\X}{\mathcal{X}}
\newcommand{\Y}{\mathcal{Y}}

\newcommand{\PX}{\mathcal{P}(\X)}
\newcommand{\PY}{\mathcal{P}(\Y)}
\newcommand{\PR}{\mathcal{P}(\R)}
\newcommand{\PRp}{\mathcal{P}(\R^+)}

\renewcommand{\H}[2]{\ensuremath{\operatorname{H}\left(\left.#1\vphantom{#2}\right|#2\vphantom{#1}\right)}}
\newcommand{\Hnm}{\operatorname{H}(\nu\mid\mu)}
\newcommand{\Hnb}{\operatorname{H}(\nu\mid\mu)+\operatorname{H}(\beta\mid\mu)}
\newcommand{\ent}{\operatorname{Ent}}
\newcommand{\Var}{\operatorname{Var}}

\renewcommand{\a}{\alpha}

\newcommand{\T}{\mathbb{T}}
\newcommand{\ST}{\mathbb{ST}}

\newcommand{\V}{\mathcal{V}}
\newcommand{\A}{\mathcal{A}}

\newcommand{\lm}{\mathcal{L}ip_{\sharp}\mu_1}
\newcommand{\ucm}{\mathcal{UC}_{\sharp}\mu_1}

\begin{document}
\title[Characterization of Talagrand's like transportation-cost inequalities \ldots]{Characterization of Talagrand's like transportation-cost inequalities on the real line.}
\author{Nathael Gozlan}

\date{\today}

\address{Laboratoire d'Analyse et de Math\'ematiques Appliqu\'ees (CNRS UMR 8050), Universit\'e de Marne-la-Vall\'ee, F-77454 Marne -la-Vall\'ee Cedex 2, France}
\email{nathael.gozlan@univ-mlv.fr}

\maketitle
\begin{center}
 \textsc{Universit\'e de Marne-la-Vall\'ee}
 \end{center}
\begin{abstract}
In this paper, we give necessary and sufficient conditions for Talagrand's like transportation cost inequalities on the real line. This brings a new wide class of examples of probability measures enjoying a dimension-free concentration of measure property.
Another byproduct is the characterization of modified Log-Sobolev inequalities for Log-concave probability measures on $\R$. 
\end{abstract}

\section{Introduction}
\subsection{Transportation-cost inequalities}
This article is devoted to the study of probability measures on the real axis satisfying some kind of transportation-cost inequalities.
These inequalities relate two quantities : on the one hand, an optimal transportation cost in the sense of Kantorovich and on the other hand, the relative entropy (also called Kullback-Leibler distance). Let us recall that if $\alpha : \R\to\R^+$ is a continuous even function, the optimal transportation-cost to transport $\nu\in \PR$ on $\mu\in \PR$ (the set of all probability measures on $\R$) is defined by :
\begin{equation}\label{opt.transp.cost}
\mathcal{T}_\a(\nu,\mu)=\inf_{\pi\in P(\nu,\mu)}\iint_{\R\times \R} \a(x-y)\,\pi(dxdy),
\end{equation}
where $P(\nu,\mu)$ is the set of all the probability measures on $\R\times \R$ such that $\pi(dx\times\R)=\nu$ and $\pi(\R\times dy)=\mu$. The relative entropy of $\nu$ with respect to $\mu$ is defined by
\begin{equation}\label{rel.entr}
\Hnm=\la\begin{array}{ll}\int \log \frac{d\nu}{d\mu}\,d\nu& \text{if } \nu\ll\mu\\+\infty&\text{otherwise}\end{array}\right.
\end{equation}
One will say that $\mu\in\PR$ satisfies the \emph{transportation-cost inequality with the cost function $(x,y)\mapsto\a(x-y)$} (TCI) if
\begin{equation}\label{TCI}
\forall \nu\in \PR,\quad \mathcal{T}_{\a}(\nu,\mu)\leq \Hnm,
\end{equation}
Transportation-cost inequalities of the form (\ref{TCI}) were introduced by K. Marton in \cite{Mar86,Mar96} and M. Talagrand in \cite{Tal96a}. After them, several authors studied inequality (\ref{TCI}), possibly in a multidimensional setting, for particular choices of the cost function $\a$ (see for example \cite{BG99}, \cite{BGL01}, \cite{CaGui}, \cite{OV00} or \cite{GGM05}).
The best known example of transportation-cost inequality is the so-called $\T_2$-inequality (also called Talagrand's inequality). It corresponds to the choice $\a(x)=\frac{1}{a}x^2$. One says that $\mu$ satisfies $\T_2$ with the constant $a$ if
\begin{equation}\label{T2}
\forall \nu\in \PR,\quad \mathcal{T}_2(\nu,\mu) \leq a\Hnm,
\end{equation}
writing $\mathcal{T}_2(\nu,\mu)$ instead of $\mathcal{T}_{x^2}(\nu,\mu)$.
\subsection{Links with the concentration of measure phenomenon}
The reason of the increasing interest to TCI is their links with the concentration of measure phenomenon. Roughly speaking, a probability measure which satisfies a TCI, also satisfies a \emph{dimension free} concentration of measure property. This link was first pointed out by K. Marton in \cite{Mar86}. For example, Talagrand's inequality is related to dimension-free gaussian concentration. If $\mu$ satisfies (\ref{T2}), then
\begin{equation}\label{gauss.conc}
\forall n\in \N^*,\forall A\subset \R^n \text{ measurable},\quad\forall r\geq r_A:=\sqrt{-\log \mu^n(A)} ,\qquad \mu^n\lp A^r\rp\geq 1-e^{-\frac{1}{a}(r-r_A)^2},
\end{equation}
where $\displaystyle{A^r=\la x \in \R^n \text{ such that } \exists y \in A \text{ with } |x-y|_2\leq r\ra}$ and $|\cdot|_2$ is the usual euclidean norm.\\

Replacing the function $x^2$ by an other convex function, it is possible to obtain different types of dimension free concentration estimates. For example, if $\mu$ is a probability measure which satisfies the transportation cost inequality
\begin{equation}\label{Tp}
\forall \nu\in \PR,\quad \mathcal{T}_{\a_p}(\nu,\mu)\leq a\Hnm,
\end{equation}
where $\a_p(x)=\la\begin{array}{ll}\min(|x|^2,|x|^p)& \text{if } p\in [1,2[\\
|x|^p& \text{if } p\geq 2\end{array} \right., \forall x\in \R$ then, it can be shown that
\begin{equation}\label{ap.conc}
\forall n\in \N^*,\forall A\subset \R^n \text{ measurable},\forall r\geq r_A:=\a_p^{-1}\lp-\log \mu^n(A)\rp ,\qquad \mu^n\lp A^r\rp\geq 1-e^{-\frac{1}{a}\alpha_p(r-r_A)},
\end{equation}
with $\displaystyle{A^r=\la x \in \R^n \text{ such that } \exists y \in A \text{ with } |x-y|_{\max(p,2)}\leq r\ra}$, denoting $|x|_p=\sqrt[p]{\sum_{i=1}^n|x_i|^p}$.
The probability measure $d\mu_p(x)=e^{-|x|^p}\frac{dx}{Z_p}, p\geq 1$ on $\R$ satisfies the TCI (\ref{Tp}).
The cases $p=1$ and $p=2$ were obtained by Talagrand in \cite{Tal96a}, the case $p\in (1,2)$ was treated by Gentil, Guillin and Miclo in \cite{GGM05} and the case $p\geq 2$ by Bobkov and Ledoux in \cite{BL00}.
\subsection{Strong transportation-cost inequalities}
When dealing with other cost functions than the $\a_p$'s, it is convenient to study a stronger form of the transportation-cost inequality (\ref{TCI}).

A probability measure $\mu$ will be said to satisfy the \emph{strong transportation-cost inequality with cost function $(x,y)\mapsto\a(x-y)$} (strong TCI) if
\begin{equation}\label{str.TCI}
\forall \nu,\beta\in \PR,\quad \mathcal{T}_{\a}(\nu,\beta)\leq \Hnb.
\end{equation}
Note that this inequality is a sort of symmetrized version of the usual TCI (\ref{TCI}). Of course, as $\H{\mu}{\mu}=0$,
\[\mu \text{ satisfies } (\ref{str.TCI})\quad\Rightarrow\quad \mu \text{ satisfies } (\ref{TCI}).\]
When $\a$ is convex, these two inequalities are equivalent up to constant factors. Namely, if $\a$ is convex one has
\[\mu \text{ satisfies } (\ref{TCI})\quad\Rightarrow\quad \mu \text{ satisfies the strong TCI with the cost function } (x,y)\mapsto 2\a\lp \frac{x-y}{2}\rp .\]
This elementary fact is proved in Proposition \ref{TCI->STCI}.

Strong TCIs are not new. The strong TCI (\ref{str.TCI}) is in fact equivalent to an infimal-convolution inequality. Infimal-convolution inequalities were introduced by B. Maurey in \cite{Mau91}. The translation of (\ref{str.TCI}) in terms of infimal-convolution inequalities will be stated in Theorem \ref{dual}.

These strong TCIs are powerful tools for deriving general dimension free concentration properties. Indeed, if $\mu$ verifies the strong TCI (\ref{str.TCI}), then
\begin{equation}\label{free.conc.meas}
\forall n\in \N^*,\forall A\subset \R^n \text{ measurable},\quad\forall r\geq0,\qquad \mu^n\lp A^r_{\a}\rp\geq 1-\frac{1}{\mu^n(A)}e^{-r},
\end{equation}
where $\displaystyle{A^r_\alpha=\la x\in \R^n : \exists y \in A \text{ such that } \sum_{i=1}^n \a(|x_i-y_i|)\leq r\ra}$.\\
Note that in (\ref{free.conc.meas}), the blow up $A^r_\a$ is not generated by a norm. When dealing with the $\a_p$'s, one can show that (\ref{free.conc.meas}) implies (\ref{ap.conc}).

The aim of this paper is to give general criteria guarantying that a probability measure satisfies a (strong) TCI. Before presenting our results let us recall some results of the literature.
\subsection{TCI and Logarithmic-Sobolev type inequalities}
The classical approach to study TCIs is to relate them to other functional inequalities such as Logarithmic-Sobolev inequalities.
The main work on the subject is the article by F. Otto and C. Villani on Talagrand's inequality (see \cite{OV00}). They proved that if $\mu\in \mathcal{P}(\R)$ satisfies the Logarithmic-Sobolev inequality
\[\ent_\mu (f^2)\leq C\int f'^2\,d\mu, \quad \forall f\]
then it satisfies Talagrand's inequality (\ref{T2}) with the same constant $C$. In fact, this result is true in a multidimensional setting. Soon after Otto and Villani, S.G. Bobkov, I. Gentil and M. Ledoux provided an other proof of this result (see \cite{BGL01}).

Different authors have tried to generalize this approach to study TCIs associated to other cost functions. Let us summarize these results. Define \[\theta_p(x)=\la\begin{array}{ll}x^2& \text{if } |x|\leq 1\\
\frac 2 p |x|^p+1-\frac{2}{p}& \text{if } |x|\geq 1\end{array} \right.,\forall p\in[1,2] 
\]
The function $\theta_p$ is just a convex function resembling to the previously defined $\a_p$. Let $\theta_p^*$ be the convex conjugate of $\theta_p$, which is defined by
\[\theta_p^*(y)=\sup_{x\in \R}\la xy - \theta_p(x)\ra.\]
If $\mu$ satisfies the following modified Logarithmic-Sobolev inequality
\begin{equation}\label{mod.Log-Sob}
\ent_\mu(f^2)\leq C\int \theta_p^*\lp\frac{tf'}{f}\rp f^2\,d\mu,\quad \forall f.
\end{equation}
for some $C,t>0$, then $\mu$ satisfies the TCI (\ref{Tp}) for some constant $a>0$.

When $p=2$, one recovers Otto and Villani's result.
The case $p=1$ was treated by Bobkov, Gentil and Ledoux in \cite{BGL01}. Note that in this case, the inequality (\ref{mod.Log-Sob}) is equivalent to Poincar\'e inequality (see \cite{BL97}).
The case $p\in (1,2)$ is due to I. Gentil, A. Guillin and L. Miclo see \cite{GGM05}.

Now the question is to know if the TCI (\ref{Tp}) is equivalent to the modified Log-Sobolev (\ref{mod.Log-Sob}). Here are some elements of answer :\\
- This is true for $p=1$. When $p=1$, inequalities (\ref{Tp}) and (\ref{mod.Log-Sob}) are both equivalent to Poincar\'e inequality (see \cite{BL97} and \cite{BGL01}).\\
- This is true as far as Log-concave distributions are concerned (see Corollary 3.1 of \cite{OV00} or Theorem 2.9 of \cite{GGM05}).\\
- For $p=2$, P. Cattiaux and A. Guillin have furnished in \cite{CaGui} an example of a probability measure which does not satisfy the Logarithmic-Sobolev inequality but satisfies Talagrand's inequality. To construct their counterexample, they give an interesting sufficient condition for Talagrand's inequality on the real line. They proved that a probability measure $\mu\in \PR$ of the form $d\mu=e^{-V}\,dx$ satisfies Talagrand's inequality (\ref{T2}) for some constant $a>0$, as soon as the potential $V$ satisfies the following condition :
\begin{equation}\label{Cat.Guil}
\limsup_{x\rightarrow \pm \infty} \frac{x}{V'(x)} <+\infty.
\end{equation}
\subsection{Presentation of the results}
In this paper, we will give necessary and sufficient conditions under which a probability measure $\mu$ on $\R$ satisfies a strong TCI. We will always assume that $\mu$ has no atom ($\mu \{x\}=0$ for all $x\in \R$) and full support ($\mu(A)>0$ for all open set $A\subset\R$).

First let us define the set of admissible cost functions. During the paper, $\mathcal{A}$ will be the class of all the functions $\a : \R\to\R^+$ such that
\begin{itemize}
\item $\a$ is even,
\item $\a$ is a continuous function, nondecreasing on $\R^+$ with $\a(0)=0$,
\item $\a$ is super-additive on $\R^+$ : $\a(x+y)\geq \a(x)+\a(y)$, $\forall x,y\geq 0$,
\item $\a$ is quadratic near $0$ : $\a(t)=|t|^2,\forall t\in [-1,1]$.
\end{itemize}
One will write $\mu\in \T_\a(a)$ (resp. $\mu\in\ST_\a(a)$) if $\mu$ satisfies the TCI (resp. the strong TCI) with the cost function $(x,y)\mapsto\a(a(x-y))$.

\subsubsection{The main result} Our main result (Theorem \ref{char.lm}) characterizes the strong TCIs on a large class $\lm\subset \PR$. Roughly speaking this set is the class of all probability measures which are Lipschitz deformation of the exponential probability measure $d\mu_1(x)=\frac{1}{2} e^{-|x|}\,dx$. More precisely $\mu$ is in $\lm$ if the monotone rearrangement map $T$ transporting $\mu_1$ on $\mu$ is Lipschitz. This map $T$ is defined by $T(x)=F^{-1}\circ F_1$, where $F$ (resp. $F_1$) is the cumulative distribution function of $\mu$ (resp. $\mu_1$), and is such that $\mu=T\sharp\mu_1$, where $T\sharp\mu_1$ denotes the image of $\mu_1$ under $T$.

If $\mu$ belongs to $\lm$ then one has the following characterization : $\mu$ satisfies the strong TCI $\ST_\a(a)$ for some constant $a>0$ if and only if there is some $b>0$ such that
\[K^+(b):=\sup_{x\geq m}\int e^{\a(bu)}\,d\mu_x^+(u)<+\infty\qquad\text{and}\qquad K^-(b):=\sup_{x\leq m}\int e^{\a(bu)}\,d\mu_x^-(u)<+\infty,\]
where $m$ is the median of $\mu$ and where $\mu_x^+$ and $\mu_x^-$ are probability measures on $\R^+$ defined as follows :
\[\mu_x^+=\mathcal{L}(X-x|X\geq x)\qquad\text{and}\qquad\mu_x^-=\mathcal{L}(x-X|X\leq x),\]
with $X$ a random variable of law $\mu$.

The result furnished by Theorem \ref{char.lm} is quite satisfactory. Firstly, though partial, this result covers all the 'regular' cases. Namely, it can be shown that if $\mu$ satisfies $\ST_a(a)$ then it satisfies a spectral gap inequality. But the elements of $\lm$ satisfy the spectral gap inequality too. Examples of probability measures not belonging to $\lm$ but satisfying the spectral gap inequality are known but are rather pathological... Secondly, one can easily derive from the above result an explicit sufficient condition for probability measures $d\mu=e^{-V}\,dx$ with $V$ satisfying a certain regularity condition.

Before stating this explicit sufficient condition, one needs to introduce the class of 'good' potentials $V$. Let $\mathcal{V}$ be the set of function $f:\R\to \R$ of class $\mathcal{C}^2$ such that
\begin{itemize}
\item there is $x_o>0$ such that $f'>0$ on $(-\infty,-x_o]\cup [x_o,+\infty),$
\item $\displaystyle{\frac{f''(x)}{f'^2(x)}\xrightarrow[x\rightarrow \pm \infty]{}0}$.
\end{itemize}
In Theorem \ref{a'/V'}, we prove that if $d\mu=e^{-V}\,dx$ with $V\in \mathcal{V}$ and $\a \in \A\cap \V$, then $\mu$ satisfies $\ST_\a(a)$ for some constant $a>0$ as soon as the following conditions hold :
\[\liminf_{x\rightarrow \pm \infty}|V'(x)|>0\qquad\text{and}\qquad \exists \lambda>0\text{ such that }\limsup_{x\to \pm \infty} \frac{\a'(\lambda x)}{V'(x+m)}<+\infty.\]
The first condition guaranties that $\mu$ belongs to $\lm$ and the second one that $K^+(b)<+\infty$ and $K^-(b)<+\infty$ for some positive $b$. This sufficient condition completely extends the result by P. Cattiaux and A. Guillin concerning $\T_2$. Our approach is completely different.

\subsubsection{The particular case of Log-concave distributions}
A particularly nice case is when $\mu$ is Log-concave. Recall that $\mu$ is said to be Log-concave if $\log (1-F)$ is concave, $F$ being the cumulative distribution function of $\mu$. If $\mu$ is Log-concave then it belongs to $\lm$. Furthermore, $\mu$ satisfies $\ST_\a(a)$ for some constant $a>0$ if and only if
\begin{equation}\label{mom}
\exists b>0,\quad \int e^{\a(bx)}\,d\mu(x)<+\infty.
\end{equation}

This result enables us to derive sufficient conditions for modified Logarithmic Sobolev inequalities. Using well known techniques, we prove in Theorem \ref{TCI->Log-Sob} that if $\mu$ is a Log-concave distribution which satisfies the inequality $\ST_\a(a)$ then $\mu$ satisfies the following modified Log-Sobolev inequality
\begin{equation}\label{mod.Log-Sob2}
\ent_\mu(f^2)\leq C\int \alpha^*\lp t\frac{f'}{f}\rp f^2\,d\mu,\quad \forall f,
\end{equation}
for some $c,t>0$.
Consequently, if the Log-concave distribution $\mu$ satisfies the moment condition (\ref{mom}), it satisfies the modified Log-Sobolev inequality (\ref{mod.Log-Sob2}) (see Theorem \ref{suf.log-sob1} and Corollary \ref{suf.log-sob2}). This extends and completes the results of Gentil, Guillin and Miclo (see \cite{GGM05} and \cite{GGM06}).
\subsubsection{A word on the method} The originality of this paper is that transportation cost inequalities are studied without the help of Logarithmic-Sobolev inequalities. Our results rely on a simple but powerful perturbation method which is explained in section 3. Roughly speaking, we show that if $\mu$ satisfies some (strong) TCI then $T_\sharp\mu$ satisfies a (strong) TCI with a skewed cost function. This principle enables us to derive new (strong) TCIs from old ones. More precisely if $\mu_\rf$ is a known probability measure satisfying some (strong) TCI and if one is able to construct a map $T$ transporting $\mu_\rf$ on an other probability measure $\mu$, then $\mu$ will satisfy a (strong) TCI too. This principle is true in any dimension. The reason why this paper deals with dimension one only is that the optimal-transportation of measures is extremely simple in this framework. 
\tableofcontents

{\bf Acknowledgements.} I want to warmly acknowledge Christian L\'eonard and Patrick Cattiaux for so many interesting 
conversations on functional inequalities and other topics.

\section{Preliminary results}
In this section, we are going to recall some well known results on TCI, namely their dual translations, their tensorization properties and their links with the concentration of measure phenomenon.

\paragraph{\bf General Framework :}
Most of the forthcoming results are available in a very general framework which we shall now describe.

Let $\X$ be a polish space and let $c:\X\times\X\to\R^+$ be a lower semi-continuous function, called the \emph{cost function}.
The set of all the probability measures on $\X$ will be denoted by $\PX$.
The optimal transportation cost between $\nu\in \PX$ and $\mu\in \PX$ is defined by
\[\mathcal{T}_c(\nu,\mu)=\inf_{\pi \in P(\nu,\mu)}\iint_{\X\times \X}c(x,y)\,\pi(dxdy),\]
where $P(\nu,\mu)$ is the set of all the probability measures on $\X\times \X$ such that $\pi(dx\times \X)=\nu$ and $\pi(\X \times dy)=\mu$.

A probability measure $\mu$ is said to satisfy the TCI with the cost function $c$ if
\begin{equation}\label{TCI.abstract}
\forall \nu\in \PX,\quad \mathcal{T}_c(\nu,\mu)\leq \Hnm.
\end{equation}
A probability measure $\mu$ is said to satisfy the strong TCI with the cost function $c$ if
\begin{equation}\label{str.TCI.abstract}
\forall \nu,\beta\in \PX,\quad \mathcal{T}_c(\nu,\beta)\leq \Hnb.
\end{equation}
\subsection{TCI vs Strong TCI}
\begin{prop}\label{TCI->STCI}
Let $\X=\R^p$ and suppose that $c(x,y)=\theta(x-y)$, with $\theta:\R^p\to\R^+$ a convex function such that $\theta(-x)=\theta(x)$.
If $\mu$ satisfies the TCI with the cost function $c$ then $\mu$ satisfies the strong TCI with the cost function $\tilde{c}$ defined by \[\tilde{c}(x,y)=2\theta\lp\frac{x-y}{2}\rp,\quad \forall x,y\in \R^p.\]
\end{prop}
\proof
Let $\pi_1\in P(\nu,\mu)$ and $\pi_2\in P(\mu,\beta)$. One can construct $X,Y,Z$ three random variables such that $\mathcal{L}(X,Y)=\pi_1$ and $\mathcal{L}(Y,Z)=\pi_2$ (see for instance the Gluing Lemma of \cite{Vill} p. 208).
Thus, using the convexity of $\theta$, one has
\begin{align*}
\mathcal{T}_{\tilde{c}}(\nu,\beta)&\leq\E\lc2\theta\lp \frac{X-Z}{2}\rp\rc\leq\E\lc\theta\lp X-Y\rp\rc +\E\lc\theta\lp Y-Z\rp\rc\\
&=\int c(x,y)\,\pi_1(dxdy)+\int c(y,z)\,\pi_2(dydz).
\end{align*}
Optimizing in $\pi_1$ and $\pi_2$ yields
\[\mathcal{T}_{\tilde{c}}(\nu,\beta)\leq \mathcal{T}_c(\nu,\mu)+\mathcal{T}_{c}(\beta,\mu),\quad\forall \nu,\beta\in\mathcal{P}(\R^p).\]
Consequently, if $\mu$ satisfies the TCI with the cost function $c$ then $\mu$ satisfies the strong TCI with the cost function $\tilde{c}$.
\endproof
Let $\theta : \R^p\to\R^+$ be a symmetric function ($\theta(-x)=\theta(x)$). One will say that a probability measure $\mu$ on $\R^p$ satisfies the inequality $\ST_\theta(a)$ if it satisfies the strong TCI with the cost function $c(x,y)=\theta(a(x-y))$.
\begin{lem}
Let $\theta$ be as above and suppose that $\theta(kx)\geq k\theta(x), \forall k\in \N, \forall x\in \R^p$. Let $b_1,b_2>0$ and define $\tilde{\theta}(x)=b_1\theta(b_2x),\forall x\in \R^p$. Then, $\mu$ satisfies $\ST_\theta(a)$ for some $a>0$ if and only if $\mu$ satisfies $\ST_{\tilde{\theta}}(\tilde{a})$ for some $\tilde{a}>0$.
\end{lem}
\proof (See also the proof of Corollary 1.3 of \cite{Tal96a})
Suppose that $\mu$ satisfies $\ST_\theta(a)$ for some $a>0$.
Let $k\in \N$ \st $k\geq b_1$, then $\theta(ax)\geq k\theta(ax/k)\geq \tilde{\theta}(ax/(b_2k))$. Hence, $\mu$ satisfies $\ST_{\tilde{\theta}}\lp\frac{a}{b_2k}\rp.$
\endproof
\subsection{Links with the concentration of measure phenomenon}
The following Theorem explains how to deduce concentration of measure estimates from a strong TCI. The argument used in the proof is due to K. Marton and M. Talagrand see (\cite{Mar86} and the proof of Corollary 1.3 of \cite{Tal96a}).
\begin{thm}\label{conc}Let $(X,d)$ be a polish space and $c:\X\times \X\to \R^+$ be a continuous cost function.
Suppose that $\mu$ satisfies the strong TCI with the cost function $c$, then
\begin{equation}\label{conc.str.TCI.abs}
\forall A\subset \X \text{ measurable},\quad\forall r\geq 0,\qquad \mu\lp A_c^r\rp\geq 1-\frac{1}{\mu(A)}e^{-r},
\end{equation}
where
$A_c^r=\la y\in \X : \exists x \in A \text{ such that } c(x,y)\leq r\ra.$
\end{thm}
\proof
Let $A, B\in \X$ and define $\mu_A=\frac{\mu(\,\cdot\,\cap A)}{\mu(A)}$ and $\mu_B=\frac{\mu(\,\cdot\,\cap B)}{\mu(B)}$. Since $\mu$ satisfies the strong TCI, one has :
\begin{equation}\label{Marton.str.TCI}
c(A,B)\leq\mathcal{T}_c(\mu_A,\mu_B)\leq \H{\mu_A}{\mu} + \H{\mu_B}{\mu} = -\log\mu(A) -\log\mu(B),
\end{equation}
with $c(A,B)=\inf\la c(x,y) : x\in A, y\in B\ra.$
Now taking $B=\X-A^r_c$ in (\ref{Marton.str.TCI}) yields the desired result.
\endproof
\subsection{Dual translation of transportation-cost inequalities}
\subsubsection{Kantorovich Rubinstein Theorem and its consequences}According to the celebrated Kantorovich-Rubinstein Theorem, optimal transportation costs admit a dual representation which is the following :
\begin{equation}\label{Kant.Rub}
\forall \nu,\mu \in \PX,\quad \mathcal{T}_c(\nu,\mu)=\sup_{(\psi, \varphi)\in \Phi_c}\la\int\psi\,d\nu-\int \varphi\,d\mu\ra,
\end{equation}
where $\Phi_c=\la (\psi,\varphi) \in B(\X)\times B(\X) : \psi(x)-\varphi(y)\leq c(x,y), \forall x,y\in \X\ra$ and $B(X)$ is the set of bounded measurable functions on $\X$. The dual representation (\ref{Kant.Rub}) is in particular true if $c$ is lower semi-continuous function defined on a polish space $\X$ (see for instance Theorem 1.3 \cite{Vill}). Furthermore, $B(\X)$ can be replaced by $C_b(\X)$, the set of bounded continuous functions on $\X$.

The infimal-convolution operator $Q_c$ is defined by \[Q_c\varphi(x)=\inf_y\la\varphi(y)+c(x,y)\ra,\] for all $\varphi \in B(\X)$. If $c$ is continuous, $x\mapsto Q_c\varphi(x)$ is measurable (in fact upper semi continuous) and it can be shown that
\begin{align*}
\forall \nu,\mu \in \PX,\quad \mathcal{T}_c(\nu,\mu)&=\sup_{\varphi \in C_b(\X)}\la\int Q_c\varphi\,d\nu-\int \varphi\,d\mu\ra,\\
&=\sup_{\varphi \in B(\X)}\la\int Q_c\varphi\,d\nu-\int \varphi\,d\mu\ra.
\end{align*}
Since optimal transportation costs admit a dual representation, it is natural to ask if TCIs and strong TCIs admit a dual translation too.
The answer is given in the following theorem.
\begin{thm}\label{dual}
Suppose that $c$ is a continuous cost-function on the Polish space $\X$.
\begin{enumerate}
\item $\mu$ satisfies the TCI (\ref{TCI.abstract}) if and only if
\begin{equation}\label{dual.TCI}
\forall \varphi\in B(\X),\quad \int e^{Q_c\varphi}\,d\mu\cdot e^{-\int \varphi\,d\mu}\leq 1.
\end{equation}
\item $\mu$ satisfies the strong TCI (\ref{str.TCI.abstract}) if and only if
\begin{equation}\label{dual.str.TCI}
\forall \varphi\in B(\X),\quad \int e^{Q_c\varphi}\,d\mu\cdot \int e^{-\varphi}\,d\mu\leq 1.
\end{equation}
\end{enumerate}
Furthermore in the preceding statements $B(\X)$ can be replaced by
$C_b(\X)$.
\end{thm}
\proof
The first point is due to S.G. Bobkov and F. G\"otze (see the proof of Theorem 1.3 and (1.7) of \cite{BG99}). The interested reader can also find an alternative proof of this result in \cite{GozLeo} (see Corollary 1). In this latter proof, Large Deviations Theory techniques are used.
One can easily adapt the one or the other approach to derive the dual version of strong TCIs (\ref{dual.str.TCI}). This is left to the reader.
\endproof
\begin{rem}As mentioned in the introduction, inequalities of the form (\ref{dual.str.TCI}) are called infimal-convolution inequalities. These inequalities were introduced by B. Maurey in \cite{Mau91}. Note that Maurey's work is anterior to the paper \cite{Tal96a} and \cite{BG99}. A good account on infimal-convolution inequalities can be found in M. Ledoux's book \cite{Led}. In this article, we have chosen to privilege the strong TCI (\ref{str.TCI.abstract}) form, which is the primal form of (\ref{dual.str.TCI}). The reason is that we find (\ref{str.TCI.abstract}) more intuitive.
\end{rem}
\subsubsection{Application : strong TCI and integrability}
Let us detail an important application of the infimal-convolution formulation of strong TCI.
\begin{prop}
Let $c$ be a continuous cost function on the Polish space $\X$. Suppose that $\mu\in \PX$ satisfies the strong TCI with the cost function $c$. Let $A\subset \X$ be a measurable set and define $c(x,A)=\inf_{y\in A}c(x,y)$. One has
\begin{equation}\label{integrability}
\int e^{c(x,A)}\,d\mu(x)\cdot\mu(A)\leq 1.
\end{equation}
\end{prop}
\begin{rem}
This integrability property was first noticed by B. Maurey in \cite{Mau91}. Note that the inequality (\ref{integrability}) implies the concentration estimate (\ref{conc.str.TCI.abs}).
\end{rem}
\proof
Define, for all $p\in \N$, $\varphi_A^p(x)=\la\begin{array}{ll} 0 & \text{if } x\in A\\
p& \text{if } x\in A^c \end{array}\right..$
As $\varphi_A^p$ is bounded, one can apply (\ref{dual.str.TCI}), this yields
\[\int e^{Q_c\varphi_A^p}\,d\mu\cdot \int e^{-\varphi_A^p}\,d\mu\leq 1.\]
An easy computation shows that $Q_c\varphi_A^p(x)=\min(c(x,A),p)\xrightarrow[p\rightarrow + \infty]{}c(x,A)$ and $e^{-\varphi_A^p}\xrightarrow[p\rightarrow + \infty]{}\1_A$. Using the monotone convergence theorem, one gets he desired inequality.
\endproof
The following Corollary will be very useful in the sequel.
\begin{cor}\label{integrability2}
Let $\mu$ be a probability measure on $\R$ satisfying the strong TCI with the cost function $c(x,y)=\alpha(x-y)$, with $\alpha$ a continuous symmetric non decreasing function. For all $x\in \R$, define
\[\mu_x^+=\mathcal{L}(X-x|X\geq x),\qquad\text{and}\qquad\mu_x^-=\mathcal{L}(x-X|X\leq x),\]
where $X$ is a random variable with law $\mu$.
Then,
\begin{align*}
&\int_0^{+\infty}e^{\a}\,d\mu_x^+\leq\frac{1}{\mu(-\infty,x]}+1,\forall x\in\R\\
&\int_0^{+\infty}e^{\a}\,d\mu_x^-\leq\frac{1}{\mu[x,+\infty)}+1,\forall x\in\R
\end{align*}
In particular,
\[\int e^{\a}\,d\mu\leq \frac{1}{\mu(\R^+)\mu(\R^-)}-1.\]
\end{cor}
\proof
Let $A=(-\infty,x]$. It is easy to show that $c(y,A)=\a(y-x)$ if $y\geq x$ and $0$ else.
Applying (\ref{integrability}) with this $A$ yields
\[\lp\mu(-\infty,x]+\int_x^{+\infty}e^{\a(y-x)}\,d\mu(y)\rp\cdot \mu(-\infty,x]\leq 1.\]
Rearranging the terms, one gets
\[\int_x^{+\infty}e^{\a(y-x)}\,d\mu(y)\leq \frac{1-\mu(-\infty,x]^2}{\mu(-\infty,x]}.\]
Dividing both sides by $\mu[x,+\infty)$ gives the result. Working with $A=[x,+\infty)$ gives the integrability property for $\mu_x^-$.
Now,
\begin{align*}
\int e^\a\,d\mu& =\mu(\R^+)\int_0^{+\infty} e^\a\,d\mu_0^++\mu(\R^-)\int_0^{+\infty} e^\a\,d\mu_0^-\\
& \leq 1+\frac{\mu(\R^+)}{\mu(\R^-)}+\frac{\mu(\R^-)}{\mu(\R^+)}\\
&=\frac{1}{\mu(\R^+)\mu(\R^-)}-1.
\end{align*}
\endproof
\subsection{Tensorization property of (strong) TCIs}
If $\mu_1$ and $\mu_2$ satisfy a (strong) TCI, does $\mu_1\otimes\mu_2$ satisfy a (strong) TCI ?
The following Theorem gives an answer to this question.
\begin{thm}\label{tens}
Let $(\X_i)_{i=1\ldots n}$ be a family of Polish spaces. Suppose that $\mu_i$ is a probability measure on $\X_i$ satisfying a (strong) TCI on $\X_i$ with a continuous cost function $c_i$ such that $c_i(x,x)=0,\forall x\in \X_i$. Then the probability measure $\mu_1\otimes\cdots\otimes\mu_n$ satisfies a (strong) TCI on $\X_1\times\cdots \times\X_n$ with the cost function $c_1\oplus\cdots\oplus c_n$ defined as follows :
\[\forall x,y \in \X_1\times\cdots\times\X_n,\quad c_1\oplus\cdots\oplus c_n(x,y)=\sum_{i=1}^nc_i(x_i,y_i).\]
\end{thm}
\proof
There are two methods to prove this tensorization property. The first one is due to K. Marton and makes use of a coupling argument (the so called \emph{Marton's coupling argument}). It is explained in several places : in Marton's original paper \cite{Mar86}, in Talagrand's paper on $\T_2$ \cite{Tal96a} or in M. Ledoux book \cite{Led} (Chapter 6). The second method uses the dual forms (\ref{dual.TCI}) and (\ref{dual.str.TCI}). This approach was originally developed by B. Maurey in \cite{Mau91} for infimal-convolution inequalities (see Lemma 1 of \cite{Mau91}). In the case of TCIs, the proof is given in great details in \cite{GozLeo} (see the proof of Theorem 5).
\endproof
\begin{rem}
Several authors have obtained non-independent tensorization results for transportation-cost inequalities and related inequalities (see \cite{Mar96}, \cite{Sam00} or \cite{DGW03}).
\end{rem}
Applying Theorem \ref{tens} together with Theorem \ref{conc}, one obtains the following Corollary :
\begin{cor}\label{dim.free}
Let $c$ be a continuous cost function on the Polish space $\X$ such that $c(x,x)=0,\forall x\in \X$. Suppose that $\mu\in \PX$ satisfies the strong TCI with the cost function $c$.
Then,
\[\forall n\in \N^*,\forall A \text{ measurable}, \forall r\geq 0, \quad\mu^n(A^r_{c})\geq 1- \frac{1}{\mu^n(A)}e^{-r},\]
where $A^r_c=\la x\in \X ^n : \exists y \in A \text{ such that } \sum_{i=1}^n c(x_i,y_i)\leq r\ra.$
\end{cor}
\section{The perturbation method for (strong) TCIs}
\subsection{The contraction principle in an abstract setting}
In the sequel, $\X$ and $\Y$ will be Polish spaces.
If $\mu$ is a probability measure on $\X$ and $T:\X\to\Y$ is a measurable map, the image of $\mu$ under $T$ will be denoted by $T_\sharp \mu$, it is the probability measure on $\Y$ defined by
$$\forall A \subset \Y \text{ measurable},\quad T_\sharp\mu(A)=\mu\left(T^{-1}(A)\right).$$
In this section, we will explain how a (strong) TCI is modified when the reference probability measure $\mu$ is replaced by the image $T_\sharp \mu$ of $\mu$ under some map $T$.
\begin{thm}\label{skewed.cost}
Let $T:\X\to\Y$ be a measurable bijection. If $\mu_\rf$ satisfies the (strong) TCI with a cost function $c_\rf$ on $\X$, then $T_\sharp\mu_\rf$ satisfies the (strong) TCI with the cost function $c_\rf^T$ defined on $\Y$ by
\[c_\rf^T(y_1,y_2)=c_\rf(T^{-1}y_1,T^{-1}y_2),\quad \forall y_1,y_2\in \Y.\]
\end{thm}
In other word, $T_\sharp\mu_\rf$ satisfies the (strong) TCI with a \emph{skewed cost function}.
\proof
Let us define $Q(y_1,y_2)=(T^{-1}y_1,T^{-1}y_2), \forall y_1,y_2\in \Y$.
Let $\nu,\beta\in \PY$ and take $\pi\in P(\nu,\beta)$, then $\int c_\rf^T(y_1,y_2)\,d\pi=\int c(x,y)\,dQ_\sharp \pi,$ so $\displaystyle{\mathcal{T}_{c_\rf^T}(\nu,\beta)=\inf_{\pi\in Q_\sharp P(\nu,\beta)}\int c(x,y)d\pi.}$
But it is easily seen that $Q_\sharp P(\nu,\beta)=P(T^{-1}_\sharp \nu,T^{-1}_\sharp\beta)$. Consequently \[\mathcal{T}_{ c_\rf^T}(\nu,\beta)=\mathcal{T}_{c_\rf}(T^{-1}_\sharp \nu,T^{-1}_\sharp\beta).\]
If $\mu_\rf$ satisfies the strong TCI with the cost function $c_\rf$, then
\[\mathcal{T}_{c_\rf}(T^{-1}_\sharp \nu,T^{-1}\beta)\leq \H{T^{-1}_\sharp \nu}{\mu_\rf}+\H{T^{-1}_\sharp \beta}{\mu_\rf}\]
But \[\H{T^{-1}_\sharp \nu}{\mu_\rf}=\H{T^{-1}_\sharp \nu}{T^{-1}_\sharp  T_\sharp  \mu_\rf}=\H{\nu}{T_\sharp  \mu_\rf},\]
where the last equality comes from the following classical invariance property of relative entropy : $\H{S_\sharp \nu_1}{S_\sharp \nu_2}=\H{\nu_1}{\nu_2}$.
Hence
\[\forall \nu,\beta \in \PY,\quad \mathcal{T}_{c_\rf^T}(\nu,\beta)\leq \H{\nu}{T_\sharp\mu_\rf}+\H{\beta}{T_\sharp\mu_\rf}.\]
\endproof
The Corollary bellow explains the method we will use in the sequel to derive new (strong) TCIs from known ones.
\begin{cor}[Contraction principle]\label{contr.princ}
Let $\mu_{\mathrm{ref}}$ be a probability measure on $\X$ satisfying a (strong) TCI with a continuous cost function $c_\rf$. In order to prove that a probability measure $\mu$ on $\Y$ satisfies the (strong) TCI with a continuous cost function $c$, it is enough to build an application $T:\X\to\Y$ such that $\mu=T_\sharp\mu_\rf$ and
\[c(Tx_1,Tx_2)\leq c_\rf(x_1,x_2),\quad \forall x_1,x_2 \in \X.\]
\end{cor}
This contraction property of strong TCIs (written in their infimal-convolution form) was first observed by B. Maurey (see Lemma 2 of \cite{Mau91}).
\proof
We assume that $\mu_\rf$ satisfies the strong TCI with the cost function $c_\rf$.
Let $\varphi:\Y\to\R$ be a bounded map. Then, for all $x_1\in \X$
\begin{align*}
Q_c\varphi(Tx_1)    &=\inf_{y\in \Y}\{\varphi(y)+c(Tx_1,y)\}\leq \inf_{x_2\in \X}\{\varphi(Tx_2)+c(Tx_1,Tx_2)\}\\
                    &\leq \inf_{x_2\in \X}\{\varphi\circ T(x_2))+c_\rf(x_1,x_2)\}=Q_{c_\rf}\lp\varphi\circ T\rp.
\end{align*}
Thus,
\begin{align*}
\int e^{Q_c\varphi}\,d\mu\cdot\int e^{-\varphi}\,d\mu & = \int e^{Q_c\varphi}\circ T\,d\mu_\rf\cdot\int e^{-\varphi\circ T}\,d\mu_\rf\\
& \leq \int e^{Q_{c_\rf}\lp\varphi\circ T\rp}\,d\mu_\rf\cdot\int e^{-\varphi\circ T}\,d\mu_\rf\\
&\leq 1
\end{align*}
where the last inequality follows from (\ref{dual.str.TCI}).
\endproof
\begin{rem}
If $T$ is invertible, the proof above can be simplified using Theorem \ref{skewed.cost}. Namely, according to Theorem \ref{skewed.cost}, $\mu$ satisfies the (strong) TCI with the cost function $c_\rf^T$. But, by hypothesis, $c\leq c_\rf^T$, so $\mu$ satisfies the (strong) TCI with the cost function $c$.
\end{rem}
\subsection{The contraction principle on the real line}
\subsubsection{Monotone rearrangement}
We are going to apply the contraction principle to probability measures on the real line. The reason why dimension one is so easy to handle is the existence of a good map $T$ which pushes forward $\mu_\rf$ on $\mu$ : the monotone rearrangement.
\begin{thm}[Monotone rearrangement]
Let $\mu_\rf$ and $\mu$ be probability measures on $\R$ and let $F_\rf$ and $F$ denote their cumulative distribution functions :
\[F_\rf(t)=\mu_\rf(-\infty,t],\quad\forall t\in \R,\qquad \text{and}\qquad F(t)=\mu(-\infty,t],\quad\forall t\in \R.\]
If $F_\rf$ and $F$ are continuous and increasing (equivalently $\mu_\rf$ and $\mu$ have no atom and full support), then the map $T=F^{-1}\circ F_\rf$ transports $\mu_\rf$ on $\mu$, that is $T_\sharp\mu_\rf=\mu$.
\end{thm}
From now on, $T$ will always be the map defined in the preceding Theorem.
\subsubsection{About the exponential distribution} The reference probability measure $\mu_\rf$ will be the symmetric exponential distribution $\mu_1$ on $\R$ :
\[d\mu_\rf(x)=d\mu_1(x):=\frac{1}{2}e^{-|x|}\,dx.\]
\begin{thm}[Maurey, Talagrand]
The exponential measure $\mu_1$ satisfies the (strong) TCI with the cost function $\frac{1}{\kappa}c_1$, for some constant $\kappa>0$,  with $c_1$ defined by\[c_1(x,y):=\min(|x-y|,|x-y|^2),\quad \forall x,y\in \R.\]
\end{thm}
\begin{rem}~
\begin{enumerate}
\item One can take $\kappa=36$.
\item B. Maurey proved the strong TCI with the sharper cost functions $c(x,y)=\tilde{\a}_1(x-y)$, where $\tilde{\a}_1(x)=\la\begin{array}{ll}\frac{1}{36} x^2&\text{if } |x|\leq 4\\\frac{2}{9}(|x|-2)&\text{if } |x|\geq 4\end{array}\right.$ (see Proposition 1 of \cite{Mau91}). One can show that $\tilde{\a}_1\geq \frac{1}{36} \a_1$.
\item M. Talagrand proved independently that $\mu_1$ satisfies the TCI with the cost functions $c_\lambda(x,y)=\gamma_\lambda(x-y)$ where $\gamma_\lambda(x)=\lp\frac{1}{\lambda}-1\rp \lp e^{-\lambda|x|}-1+\lambda|x|\rp$ for all $\lambda\in (0,1)$ (see Theorem 1.2 of \cite{Tal96a}).
\end{enumerate}
\end{rem}
Transportation-cost inequalities associated to the cost function $c_1$ were fully characterized by I. Gentil, M. Ledoux and S . Bobkov in \cite{BGL01} in terms of Poincar\'e inequalities :
\begin{thm}[Bobkov-Gentil-Ledoux]\label{BGL}
A probability measure $\mu$ on $\R^p$ satisfies the TCI with the cost function $(x,y)\mapsto \lambda \min(|x-y|_2,|x-y|_2^2)$, for some $\lambda>0$ if and only if it satisfies a Poincar\'e inequality, that is if there is some constant $C>0$ such that
\[\Var_\mu(f)\leq C \int_{\R^p} |\nabla f|^2_2\,d\mu,\quad \forall f\]
\end{thm}
\subsubsection{Application of the contraction principle on the real line}
A good thong with the exponential distribution is that its cumulative distribution function can be explicitly computed
\begin{equation}\label{cdf.exp}
F_1(x)=\la\begin{array}{ll}1-\frac 1 2 e^{-|x|} & \text{if } x\geq 0\\
\frac 1 2 e^{-|x|} & \text{if } x\leq 0
 \end{array}\right.\quad \text{ and }\quad  F_1^{-1}(t)=\la\begin{array}{ll}-\log(2(1-t)) & \text{if } t\geq \frac 1 2\\
\log(2t) & \text{if } t\leq \frac 1 2
 \end{array}\right.
 \end{equation}
Suppose that $\mu$ is a probability measure on $\R$ having no atom and full support, then its cumulative distribution function $F$ is invertible, and the map $T$ transporting $\mu_1$ on $\mu$ can be expressed as follows :
\begin{equation}\label{T}
T(x)=\la\begin{array}{ll}F^{-1}\lp1-\frac 1 2 e^{-|x|}\rp & \text{if } x\geq 0\\
F^{-1}\lp\frac 1 2 e^{-|x|}\rp & \text{if } x\leq 0
 \end{array}\right.\quad \text{ and }\quad  T^{-1}(x)=\la\begin{array}{ll}-\log(2(1-F(x))) & \text{if } x\geq m\\
\log(2F(x)) & \text{if } x\leq m
 \end{array}\right.,
 \end{equation}
where $m$ denotes the median of $\mu$.\\

Let us introduce the following quantity :
\[\omega_{\mu}(h)=\inf\la|T^{-1}x-T^{-1}y| : |x-y|\geq h\ra,\quad \forall h\geq 0.\]

\begin{prop}
If $\mu\in \PR$ is a probability measure with no atom and full support, then $\mu$ satisfies the strong TCI with the cost function $c_\mu(x,y)=\frac 1 \kappa\a_1\circ \omega_\mu(|x-y|)$, where $\a_1(t)=\min(t,t^2),\forall t\geq 0$.
\end{prop}
\proof
By definition of $\omega_\mu$,
\[|T^{-1}x-T^{-1}y|\geq \omega_\mu(|x-y|),\quad \forall x,y\in \R.\]
Thus,
\[c_1^T(x,y)\geq \frac{1}{\kappa}\a_1\lp\omega_\mu(|x-y|)\rp,\quad \forall x,y\in \R,\]
and this achieves the proof.
\endproof
To better understand $\omega_\mu$ it is good to relate it to the continuity modulus of $T$.
\begin{defi}[The class $\ucm$]
The set of all probability measures $\mu$ on $\R$, with no atom and full support, such that the monotone rearrangement map transporting the exponential measure $\frac{1}{2}d\mu_1(x)=e^{-|x|}\,dx$ on $\mu$ is uniformly continuous is denoted by $\ucm$.
\end{defi}
The proof of the following proposition is left to the reader.
\begin{prop}
Suppose $\mu\in \ucm$, then the continuity modulus $\Delta_\mu$ of $T$ is defined by $\Delta_\mu(h)=\sup\la|Tx-Ty| : |x-y|\leq h\ra,$ $\forall h\geq 0$. It is a continuous increasing function and
\[\omega_\mu=\Delta_\mu^{-1}.\]
\end{prop}
\begin{rem}~
\begin{enumerate}
\item All the elements of $\ucm$ enjoy a dimension free concentration of measure property. Namely, if $\mu \in \ucm$, then $\mu$ satisfies the strong TCI with the cost function $c_\mu(x,y)=\a_\mu(x-y)$, where $\a_\mu(x)=\frac{1}{\kappa}\a_1\circ \omega_\mu(|x|)$. Thus
according to Corollary \ref{dim.free}, one has
\[\forall n\in \N^*, \forall A \subset \R^n,\forall r\geq0,\quad \mu^n (A^r_{\a_\mu})\geq 1- \frac{1}{\mu(A)}e^{-r},\]
 with $A^r_{\a_\mu}=\la x \in \R^n : \exists y\in A \st \sum_{i=1}^n \a_\mu(x_i-y_i)\leq r\ra.$
\item The class of all the probability measures on $\R$ satisfying a dimension free concentration of measure property is not yet identified. In \cite{BH00}, S.G. Bobkov and C. Houdr\'e studied probability measures enjoying a weak dimension free concentration property (roughly speaking one can estimate $\mu^n(A^r_\infty)$ independently of the dimension, where $A^r_\infty$ denotes the blow-up of $A$ with respect to the norm $|x|_\infty=\max_i |x_i|$).
They proved that a probability measure has this weak property if and only if the map $T$ generate a finite modulus, which means that $\Delta_\mu(h)<+\infty$ for some (equivalently for all) $h \in \R$.
\end{enumerate}
\end{rem}
In order to obtain explicit concentration properties, one has to estimate $\omega_\mu$.
\begin{prop}
Define
\begin{align*}
\omega_\mu^+(h)&=\inf\la|T^{-1}x-T^{-1}y| : |x-y|\geq h, x,y \geq m\ra\\
\omega_\mu^-(h)&=\inf\la|T^{-1}x-T^{-1}y| : |x-y|\geq h, x,y \leq m\ra
\end{align*}
then \[\omega_\mu(h)\geq \min\lp\omega_\mu^+\lp\frac h 2\rp, \omega_\mu^-\lp\frac h 2\rp\rp\]
\end{prop}
\proof
Let $x,y \in \R$ with $x\leq m\leq y$ and $y-x\geq h\geq 0$. One has
\begin{align*}
|T^{-1}y-T^{-1}x|&=T^{-1}y-T^{-1}x=\lp T^{-1}y-T^{-1}m\rp+\lp T^{-1}m-T^{-1}x\rp\\
&\geq \omega_\mu^+(y-m)+\omega_\mu^-(m-x).
\end{align*}
Since $y-x\geq h$ and $m\in [x,y]$, one has $y-m\geq \frac h 2$ or $m-x \geq \frac h 2$, thus
\[|T^{-1}y-T^{-1}x|\geq \min\lp\omega_\mu^+\lp\frac h 2\rp, \omega_\mu^-\lp\frac h 2\rp\rp.\]
\endproof
Let $X$ be a random variable with law $\mu$ and define
\begin{align}
\mu_x^+&=\mathcal{L}(X-x\mid X\geq x) \in \PRp,\quad \forall x\geq m\\
\mu_x^-&=\mathcal{L}(x-X\mid X\leq x) \in \PRp,\quad \forall x\leq m
\end{align}
In the following Proposition, the quantities $\omega_\mu^+$ and $\omega_\mu^-$ are expressed in terms of the cumulative distribution functions of the probability measures $\mu_x^+$ and $\mu_x^-$.
\begin{prop}
\[
\begin{array}{l}
\omega_\mu^+(h)=\inf\la-\log \mu_x^+[h,+\infty) : x\geq m\ra\\
\omega_\mu^-(h)=\inf\la-\log \mu_x^-[h,+\infty) : x\leq m\ra\\
\end{array},\quad \forall h\geq 0.
\]
\end{prop}
\proof
It is easy to see that $\omega_\mu^+(h)=\inf\la T^{-1}(x+h)-T^{-1}x : x\geq m\ra.$
Using (\ref{T}) one sees that \[T^{-1}(x+h)-T^{-1}x=-\log \lp \frac{1-F(x+h)}{1-F(x)}\rp=-\log \mu_x^+[h,+\infty),\]
which gives the result.
\endproof
The proof of the following Corollary is immediate.
\begin{cor}\label{unif.stoc.dom}
Let $\omega : \R^+\to \R^+$ a continuous non decreasing function with $\omega(0)=0.$ In order to show that $c_\mu(x,y)\geq \frac 1 \kappa\a_1\circ \omega\lp\frac{|x-y|}{2}\rp,$ it is enough to show that
\begin{align}
&\sup_{x\geq m}\mu_x^+[h,+\infty)\leq e^{-\omega(h)},\quad\forall h\geq 0\label{unif.stoc.dom.1}\\
&\sup_{x\leq m}\mu_x^-[h,+\infty)\leq e^{-\omega(h)},\quad\forall h\geq 0\label{unif.stoc.dom.2}
\end{align}
\end{cor}
\begin{rem}
The evolution of $\mu_x^+$ and $\mu_x^-$ with $x$ reflects the aging properties of $\mu$. Objects of this type appear naturally in reliability theory. Suppose that $X$ is nonnegative and think of $X$ as the failure time of some engine, then $\mu_x^+$ is the law of the failure after time $x$ knowing that the engine works properly at time $x$.
\end{rem}
Recall that one says that a probability measure $\nu_1$ is stochastically dominated by an other probability measure $\nu_2$ if
\[\nu_1[h,+\infty)\leq \nu_2[h,+\infty),\quad \forall h\geq 0.\]
In the sequel, this will be written $\nu_1\sd\nu_2$.
If one think of $\mu_1$ and $\mu_2$ as failure time laws, then $\nu_1 \sd \nu_2$ means that the material modeled by $\nu_2$ is more reliable than the material modeled by $\nu_1$.

It is well known that the following propositions are equivalent :
\begin{enumerate}
\item $\nu_1\sd\nu_2$,
\item $\int f\,d\nu_1\leq \int f\,d\nu_2$, for all nondecreasing $f:\R\to\R$,
\item There are $X_1$, and $X_2$ two random variables defined on the same probability space, such that $\mathcal{L}(X_1)=\nu_1$, $\mathcal{L}(X_2)=\nu_2$ and $X_1\leq X_2$ almost surely.
\end{enumerate}
Since every continuous nondecreasing function $F$ with $F(0)=0$ and $\lim_{x\rightarrow + \infty}F(x)=1$ is the cumulative distribution function of some probability measure on  $\R^+$ with no atom,
finding a function $\omega$ such that (\ref{unif.stoc.dom.1}) and (\ref{unif.stoc.dom.2}) hold is the same as finding some uniform upper bound of the probability measures $\mu_x^+$ and $\mu_x^-$ in the sense of stochastic ordering. The preceding Corollary can thus be restated as follows :
\begin{cor}\label{stoc.dom}
Let $\mu\in \PR$ be a probability measure with no atom and full support. If there is a probability measure $\mu_0\in \PRp$ with no atom such that
\[\mu_x^+\sd\mu_0,\quad\forall x\geq m\qquad\text{and}\qquad\mu_x^-\sd\mu_0,\quad\forall x\leq m,\]
then $\mu$ satisfies the strong TCI with the cost function $c$ defined by
\[c(x,y)=\frac 1 \kappa\a_1\lp-\log\lp1-F_0\lp|x-y|/2\rp\rp\rp,\quad \forall x,y\in \R,\]
where $F_0$ denotes the cumulative distribution function of $\mu_0$.
\end{cor}

\section{(Strong) TCI for Log-concave distributions}
Let $\mu\in \PR$, let $F$ be its cumulative distribution function and define $\overline{F}=1-F$. The probability measure $\mu$ is said to be \emph{Log-concave} if $\log \overline{F}$ is concave.
\subsection{A natural cost function}
Log-concave are examples of NBU (New Better than Used) distributions. This is explained in the following proposition :
\begin{prop}\label{stoc.dom.log.conc}
If $\mu\in \PR$ is a Log-concave distribution, then \[\mu_x^+\sd \mu_m^+,\quad \forall x\geq m \qquad\text{and}\qquad \mu_x^-\sd \mu_m^-,\quad \forall x\leq m.\]
\end{prop}
\proof
Let us show that $\mu_x^+\sd\mu_m^+$ for all $x\geq m$. By definition, this means that $\mu_x^+[h,+\infty)\leq \mu_m^+[h,+\infty),\forall h\geq 0,\forall x\geq m$ and this is equivalent to
\[\frac{1-F(x+h)}{1-F(x)}\leq \frac{1-F(m+h)}{1-F(m)},\quad \forall h\geq 0,\forall x\geq m.\]
Defining $\overline{F}_m^+(h)=\frac{1-F(m+h)}{1-F(m)},\forall h\geq 0$, the preceding inequality is equivalent to :
\[\overline{F}_m^+(x-m+h)\leq \overline{F}_m^+(x-m)\times \overline{F}_m^+(h),\quad \forall h\geq 0,\forall x\geq m.\]
In other word, $\mu_x^+\sd\mu_m^+$ if and only if the function $\log \overline{F}_m^+$ is sub-additive. Since $\mu$ is Log-concave, the function $\log \overline{F}_m^+$ is concave. It is easy to check that every concave function defined on $\R^+$ and vanishing at $0$ is sub-additive. This achieves the proof.
\endproof

\begin{cor}\label{natural}
If $\mu\in \PR$ is Log-concave, then it satisfies the strong TCI with the cost function $c$ defined by
\[c(x,y)=\frac 1 \kappa\a_1\lp-\log\lp G_0\lp|x-y|/2\rp\rp\rp,\quad \forall x,y\in \R,\]
with $G_0(h)=2\max \lp \overline{F}(h+m), F(-h+m)\rp,\forall h\geq 0$.\\
Furthermore, if $\mu$ is symmetric, then $\mu$ satisfies the strong TCI with the cost function
\[c(x,y)=\frac 1 \kappa\a_1\lp-\log 2\overline{F}(|x-y|/2)\rp,\quad \forall x,y\in \R.\]
\end{cor}
\proof
Let $\mu_0$ be the probability measure on $\R^+$ with cumulative distribution function $F_0=1-G_0$. According to Proposition \ref{stoc.dom.log.conc}, $\mu_x^+\sd\mu_0, \forall x\geq m$ and $\mu_x^- \sd \mu_0,\forall x\leq m$. Thus according to Corollary \ref{stoc.dom}, $\mu$ satisfies the strong TCI with the cost function \[c(x,y)=\frac 1 \kappa\a_1\lp-\log\lp G_0\lp|x-y|/2\rp\rp\rp,\quad \forall x,y\in \R.\]
Now, if $\mu$ is symmetric, then $m=0$ and $1-F(x)=F(-x),\forall x\in \R$. Consequently, $G_0(h)=2(1-F(h))$ and the result follows.
\endproof
\subsection{Characterization of (strong) TCI for Log-concave measures}
In the sequel, $\mathcal{A}$ will be the class of all the function $\a : \R\to\R^+$ such that
\begin{itemize}
\item $\a$ is even,
\item $\a$ is continuous, nondecreasing on $\R^+$ and $\a(0)=0$,
\item $\a$ is super-additive on $\R^+$ : $\a(x+y)\geq \a(x)+\a(y), \forall x,y\geq 0$,
\item $\a$ is quadratic near $0$ : $\a(t)=t^2,\forall t\in [-1,1]$.
\end{itemize}
One will say that $\mu$ satisfies the inequality $\T_\a(a)$ (resp. $\ST_\a(a)$) if $\mu$ satisfies the TCI (resp. the strong TCI) with the cost function $c(x,y)=\a(a(x-y)),\forall x,y\in \R$.
\begin{thm}\label{char.log-conc}
Let $\a\in \mathcal{A}$ and $\mu\in \PR$ a Log-concave distribution. The following propositions are equivalent
\begin{enumerate}
\item There is some constant $a>0$ such that $\mu$ satisfies the inequality $\ST_\a(a)$.
\item There is some constant $b>0$ \st $\int e^{\a(bx)}\,d\mu(x)<+\infty.$
\end{enumerate}
If $\a\in \mathcal{A}$ is convex then the same is true for TCI.
\end{thm}
\proof~\\
$[(1)\Rightarrow (2)]$ If $\mu$ satisfies $\ST_\a(a)$, then according to Corollary \ref{integrability2}, one has
\[\int e^{\a(az)}\,d\mu(z)<+\infty.\]
Hence, (2) holds with $b=a$.

$[(2)\Rightarrow (1)]$ According to Corollary \ref{natural}, $\mu$ satisfies the strong TCI with the cost function $c(x,y)=\frac 1 \kappa \a_1\lp-\log\lp G_0\lp|x-y|/2\rp\rp\rp,$ where $G_0(h)=2\max \lp \overline{F}(h+m), F(-h+m)\rp,\forall h\geq 0$.
If there is some $a$ such that
\begin{equation}\label{char.log-conc.3}
\a_1\lp-\log G_0(|x|)\rp\geq \a(ax),\forall x\in \R,
\end{equation}
then $\mu$ satisfies the strong TCI with the cost function $\frac 1 \kappa \a(a|x-y|/2)\geq \a(a|x-y|/(2\kappa))$ (since $\kappa>1$).
Hence it is enough to prove (\ref{char.log-conc.3}).
This latter condition is equivalent to
\begin{align}
2\overline{F}(m+h)=\mu_m^+[h,+\infty)& \leq e^{-\a_1^{-1}\circ \a(ah)},\quad\forall h\geq 0\label{char.log-conc.1}\\
2F(m-h)=\mu_m^-[h,+\infty)&\leq e^{-\a_1^{-1}\circ \a(ah)},\quad\forall h\geq 0.\label{char.log-conc.2}
\end{align}
We will focus on the condition (\ref{char.log-conc.1}), the same proof will work for (\ref{char.log-conc.2}). Inequality (\ref{char.log-conc.1}) is equivalent to
\[\la\begin{array}{lll} (i) & \mu_m^+[h,+\infty)\leq e^{-ah},&\quad \forall h\leq \frac 1 a\\
(ii)&\mu_m^+[h,+\infty)\leq e^{-\a(ah)},&\quad \forall h\geq \frac 1 a
\end{array}\right.\]
Let us prove that (i) holds for some $a_0>0$ and \emph{all} $h\geq 0$. Let $\varphi=\log \overline{F}$. The function $\varphi$ is concave, so
\[\varphi(m+h)\leq \varphi(m) + \varphi'_r(m)h,\quad \forall h\geq 0,\]
where $\varphi_r'(m)$ is the right derivative of $\varphi$ at point $m$. If $\varphi_r'(m)<0$, then (i) holds with $a_0=-\varphi_r'(m)$ and for all $h\geq 0$. Since $\varphi$ is increasing, $\varphi_r'(m)\leq0$. The function $\varphi$ being concave, $\varphi_r'$ is non-increasing. Consequently, if $\varphi_r'(m)=0$, then $\varphi_r'(x)=0$, for all $x\leq m$. This would imply that $\overline{F}$ is constant on $(-\infty,m]$, which is absurd.

It is clear that one can find a constant $b_0$ such that
$\int e^{\a(b_0z)}\,d\mu_m^+(z)<+\infty.$ We end the proof applying the following technical result to $\mu_m^+$.
\endproof
\begin{lem}\label{tech.lem}
Let $\nu$ be a probability measure on $\R^+$ such that
\[\nu[h,+\infty)\leq e^{-a_0h},\quad \forall h\geq 0,\]
for some $a_0>0$. Let $\a\in \mathcal{A}$ and suppose that \[\int e^{\a(b_0z)}\,d\nu(z)\leq K,\]
for some $b_0>0$ and $K>0$. Then, there is a constant $a>0$ depending only on $a_0, b_0$ and $K$ such that
\[\la\begin{array}{lll} (i) & \nu[h,+\infty)\leq e^{-ah},&\quad \forall h\leq \frac 1 a\\
(ii)&\nu[h,+\infty)\leq e^{-\a(ah)},&\quad \forall h\geq \frac 1 a
\end{array}\right.\]
\end{lem}
\proof
Using Markov's inequality, one gets
\[K\geq2e^{\a(b_0h)}\nu[h,+\infty),\quad \forall h\geq 0.\]
Thus, using the super-additivity of $\a$, one has
\[
\nu[h,+\infty)\leq Ke^{-\a(b_0h)}
 \leq \lc Ke^{-\a(b_0h/2)}\rc e^{-\a(b_0h/2)}
 \leq e^{-\a(b_0h/2)},
\]
as soon as $h\geq \frac{2}{b_0}\a^{-1}\lp\log K\rp$.
Since $\nu[h,+\infty)\leq e^{-a_0h},\quad \forall h\geq 0,$ it is now easy to check that (i) and (ii) hold with $a=\min(a_0,b_0/2,\lc\frac{2}{b_0}\a^{-1}\lp\log K\rp\rc^{-1})$.
\endproof
\subsection{Links with modified Log-Sobolev inequalities}
Recall the definition of the entropy functional :
\[\ent_\mu(f):=\int f\log f\,d\mu-\int f\,d\mu\log\int f\,d\mu.\]

\begin{defi}
Let $\beta : \R\to\R^+$ be an even convex function with $\beta(0)=0$. One says that $\mu\in \PR$ satisfies the modified Logarithmic-Sobolev inequality $LSI_\beta(C,t)$ if
\[\ent_\mu(f^2)\leq C\int \beta\lp t\frac{f'}{f}\rp f^2\,d\mu,\]
for all $f\in\mathcal{C}_c^1$ (the set of continuously differentiable functions having compact support).
\end{defi}
Note that if $\beta(x)=x^2$, one recovers the classical Logarithmic-Sobolev inequality.
The links between transportation cost inequalities and Logarithmic-Sobolev inequalities have been studied by several authors (see the works by Otto and Villani \cite{OV00}, Bobkov, Gentil and Ledoux \cite{BGL01} and more recently Gentil, Guillin and Miclo \cite{GGM05}).
The usual point of view is to prove TCI using Log-Sobolev type inequalities. Here we will do the opposite and derive Log-Sobolev inequalities from TCIs. To this end we will use the following result.
\begin{thm}\label{TCI->Log-Sob}
Let $\a \in \mathcal{A}$ be a convex function.
If $\mu=e^{-V}\,dx\in \PR$ with $V:\R\to\R$ a convex function satisfies the inequality $\T_\a(a)$, then it satisfies $LSI_{\a^*}\lp\frac{\lambda}{1-\lambda}, \frac{1}{a\lambda}\rp,$ for all $\lambda\in(0,1)$, where $\a^*$ is the convex conjugate of $\a$ :
\[\a^*(s)=\sup_{t\in \R}\la st-\a(t)\ra,\quad \forall s\in \R.\]
\end{thm}
\proof
The proof of Theorem \ref{TCI->Log-Sob} can be easily adapted from the one of Theorem 2.9 of \cite{GGM05}. The regularity issue mentioned by the authors during the proof, is irrelevant in our framework. Namely, in dimension one, the Brenier map is simply the monotone rearrangement map, and the regularity of this latter can be easily checked by hand.
\endproof
The following result follows immediately from Theorems \ref{char.log-conc} and \ref{TCI->Log-Sob}.

\begin{thm}\label{suf.log-sob1}
Let $\a \in \mathcal{A}$ be a convex function and $\mu=e^{-V}\,dx\in \PR$ with $V$ convex.
If $\int e^{\a(b|x|)}\,d\mu(x)<+\infty$, for some $b>0$, then $\mu$ satisfies the inequality $LSI_{\a^*}(C,t)$ for some $C,t>0$.
\end{thm}
\begin{rem}
Recall that the function $\theta_p$ is defined by
\[\theta_p(x)=\la\begin{array}{ll}x^2& \text{if } |x|\leq 1\\
\frac 2 p |x|^p+1-\frac{2}{p}& \text{if } |x|\geq 1\end{array} \right.,\forall p\in[1,2]
\]
In \cite{GGM05}, I. Gentil, A. Guillin and L. Miclo proved that the measure $d\mu_p(x)=\frac{1}{Z_p}e^{-|x|^p}\,dx$ with $p\in [1,2]$ satisfies the inequality $LSI_{\theta_p^*}(C,t)$ for some $C,t>0$.\\
Using classical tools, one can show that
\[\exists C,t>0 \text{ s.t. }\mu \text{ satisfies } LSI_{\theta^*_p}(C,t)\quad\Rightarrow\quad \exists b>0 \text{ s.t. }\int e^{\theta_p(ax)}\,d\mu(x)<+\infty.\]
Consequently, a Log-concave measure $\mu$ satisfies the inequality $LSI_{\theta_p^*}(C,t)$ if and only if there is some $b>0$ such that $\int e^{\theta_p(ax)}\,d\mu(x)<+\infty$.
\end{rem}
Suppose that $d\mu=e^{-V}\,dx$ with $V$ a convex and symmetric function. It is tempting to take $\a=V$ in the above theorem. To do this one has to modify the potential $V$ near $0$.
Define \[\widetilde{V}(x)=\la\begin{array}{ll} x^2 & \text{if } |x|\leq 1\\ V(a_0x) + 1-V(a_0) & \text{if } |x|\geq 1\end{array}\right.\]
Choosing $a_0>0$ such that $a_0V'(a_0)=2$ (which is always possible), one obtains a convex function. Furthermore, it is clear that one can find some $b>0$ such that $\int e^{\widetilde{V}(bx)}\,d\mu(x)<+\infty$. Applying the above Theorem, one obtains the following result
\begin{cor}\label{suf.log-sob2} With the above notations, $\mu$ satisfies the inequality $LSI_{\widetilde{V}^*}(C,t)$ for some $C,t>0$.
\end{cor}
\begin{rem}
In \cite{GGM06}, Gentil, Guillin and Miclo have obtained the preceding Corollary under the following additional assumption on $V$ :
\[\exists \ep \in [0,\frac{1}{2}],\exists M>0, \forall x\geq M,\quad (1+\ep)V(x)\leq xV'(x)\leq (2-\ep)V(x).\]
This hypothesis seems to be useless.
\end{rem}
\section{Characterization of strong TCI on a larger class of probabilities}
In this section, we give a characterization of strong TCI for probability measures belonging to a certain class $\lm$ which we shall now define.
\subsection{The Lipschitz images of the exponential measure}
\begin{defi}[The class $\lm$]
The set of all probability measures $\mu$ on $\R$, with no atom and full support, such that the monotone rearrangement map transporting the exponential measure $\frac{1}{2}d\mu_1(x)=e^{-|x|}\,dx$ on $\mu$ is Lipschitz is denoted by $\lm$.
\end{defi}
The following Proposition describes the elements of $\lm$.
\begin{prop}
Let $\mu\in \PR$ with no atom and full support and let $T$ be the monotone rearrangement map between $\mu_1$ and $\mu$.
For all $a>0$, let $\nu_a$ be the one sided exponential distribution with parameter $a$, that is $d\nu_a(z)=ae^{-ay}\1_{[0,+\infty)}(y)\,dy.$

The following assertions are equivalent
\begin{enumerate}
\item The map $T$ is $\frac 1 a$-Lipschitz.
\item The probability measures $(\mu_x^+)_{x\geq m}$ and $(\mu_x^-)_{x\leq m}$ are stochastically dominated by the exponential measure $\nu_a$ :
$\mu_x^+\leq \nu_a,\forall x\geq m$ and $\mu_x^-\sd\nu_a,\forall x\leq m.$
In other word, one has
\begin{align}
\sup_{x\geq m} \mu_x^+[h,+\infty)\leq e^{-ah},&\quad \forall h\geq 0\label{stoc.dom.exp1}\\
\sup_{x\leq m} \mu_x^-[h,+\infty)\leq e^{-ah},&\quad \forall h\geq 0\label{stoc.dom.exp2}
\end{align}
\end{enumerate}
If $\mu$ is of the form $d\mu(z)=e^{-V(z)}\,dz$ where $V$ is a continuous function, then $T$ is $\frac 1 a$-Lipschitz if and only if
\begin{equation}\label{A+A-}
A^+:=\sup_{x\geq m}\overline{F}(x)e^{V(x)}\leq \frac{1}{a} \qquad\text{and}\qquad A^-:=\sup_{x\leq m}F(x)e^{V(x)}\leq \frac{1}{a}
\end{equation}
Furthermore, if $V$ is of class $\mathcal{C}^1$, a sufficient condition for $A^+$ and $A^-$ to be finite is that
\begin{equation}
\liminf_{x\rightarrow + \infty} V'>0\qquad\text{and}\qquad \limsup_{x\rightarrow - \infty} V'<0.
\end{equation}
\end{prop}
\proof
It is easy to see that the map $T$ is $\frac{1}{a}$-Lipschitz if and only if
\begin{equation}\label{TLip}
T^{-1}z-T^{-1}y\geq a(z-y),\quad \forall z\geq y.
\end{equation}
This is equivalent to
\[T^{-1}(x+h)-T^{-1}x\geq ah,\quad \forall x\geq m,\forall h\geq 0\quad \text{ and } \quad T^{-1}x-T^{-1}(x-h)\geq ah,\quad \forall x\leq m,\forall h\geq 0.\]
Using the fact that
$T^{-1}(z)=\la\begin{array}{ll}-\log(2(1-F(z))) & \text{if } z\geq m\\
\log(2F(z)) & \text{if } z\leq m
 \end{array}\right.
$, one sees immediately that these conditions are equivalent to (\ref{stoc.dom.exp1}) and (\ref{stoc.dom.exp2}).

If $d\mu(z)=e^{-V(z)}\,dz$ with a continuous $V$, $T^{-1}$ is differentiable. Observe that (\ref{TLip}) means that $z\mapsto T^{-1}z-az$ is nondecreasing and this is equivalent to $\sup_{z\in\R }\frac{dT^{-1}}{dz}(z)\leq a$. Computing $\frac{dT^{-1}}{dz}$, one obtains immediately (\ref{A+A-}).

Finally, let us show that the condition $\displaystyle{\liminf_{x\rightarrow + \infty} V'>0}$ implies that $A^+$ is finite. Under this assumption, there is $v_0>0$ and $z_0>m$ such that for all $z\geq z_0$, one has $V'(z)\geq v_0$. If $z\geq z_0$, one thus has
$$e^{-V(z)}=\int_z^{+\infty}V'(y)e^{-V(y)}\,dy\geq v_0\int_z^{+\infty}e^{-V(y)}\,dy=v_0(1-F(y)).$$
So, $\displaystyle{\sup_{z\geq z_0}(1-F(z))e^{V(z)}\leq \frac{1}{v_0}}$. Since $\displaystyle{\sup_{m\leq z\leq z_0}(1-F(z))e^{V(z)}<+\infty}$, one concludes that $A^+<+\infty$. The same reasoning shows that the condition $\displaystyle{\limsup_{z\rightarrow - \infty} V' <0}$ implies $A^-<+\infty$.
\endproof
\begin{rem}[$\lm$ and the Poincar\'e inequality]~
\begin{enumerate}
\item According to Corollary \ref{contr.princ}, one concludes that a sufficient condition for a probability measure to satisfy the inequality $\ST_{\a_1}(a)$ for some constant $a$ is that $\mu$ belongs to $\lm$.
\item According to Theorem \ref{BGL}, a probability measure satisfies the inequality $\ST_{\a_1}(a)$ for some constant $a>0$ if and only if it satisfies the Poincar\'e inequality :
\begin{equation}\label{Poincaré}
\Var_\mu(f)\leq C\int (f')^2\,d\mu,\quad \forall f,
\end{equation}
for some constant $C>0$. Examples of probability measures not belonging to $\lm$ but satisfying a Poincar\'e inequality are known. Thus our perturbation method failed to completely characterize the inequality $\ST_{\a_1}$. Nevertheless, the above mentioned counterexamples are rather pathological, and for a large class of probability measures, $\mu$ satisfies Poincar\'e if and only if $\mu\in \lm$. This explained in the next proposition.
\end{enumerate}
\end{rem}
\begin{defi}[Good potentials]
The class $\mathcal{V}$ will be the set of all the functions $f:\R\to\R$ of class $\mathcal{C}^2$ such that
\begin{enumerate}
\item there is $x_o>0$ such that $f'>0$ on $(-\infty,-x_o]\cup [x_o,+\infty),$
\item $\displaystyle{\frac{f''(x)}{f'^2(x)}\xrightarrow[x\rightarrow \pm \infty]{}0}$.
\end{enumerate}
\end{defi}
\begin{prop}
Let $d\mu=e^{-V}\,dx$ with $V\in \V$, then
\begin{align*}
\mu \text{ satisfies Poincar\'e}&\Leftrightarrow \liminf_{x\rightarrow +\infty}V'(x)>0 \text{ and } \limsup_{x\rightarrow - \infty} V'(x)<0\\
& \Leftrightarrow\mu \in \lm.
\end{align*}
\end{prop}
\proof
According to the celebrated Muckenhoupt criterion, a probability measure $d\mu=e^{-V}dx$ with a continuous $V$ satisfies (\ref{Poincaré}) for some constant $C>0$ if and only if
\begin{equation}\label{Muckenhoupt}
D^+:=\sup_{x\geq m} \overline{F}(x)\cdot \int_m^x e^V(y)\,dy<+\infty\quad\text{and}\quad D^-:=\sup_{x\leq m} F(x)\cdot \int_x^m  e^V(y)\,dy<+\infty.
\end{equation}
Applying Proposition \ref{int.equiv}, one shows that \[\overline{F}(x)\cdot  \int_m^x e^V(y)\,dy\sim_{x\to +\infty} \frac{1}{V'^2(x)}\quad\text{and}\quad \overline{F}(x)\cdot e^{V(x)}\sim_{x\to+\infty} \frac{1}{V'(x)}.\]
From this one easily conclude that $A^+$ and $D^+$ are finite if and only if $\displaystyle{\liminf_{x\rightarrow +\infty} V'(x)>0}$.
\endproof
\subsection{Characterization of strong TCI on $\lm$}
\begin{thm}\label{char.lm}
Let $\mu\in \lm$ and $\a\in \mathcal{A}$. The following assertions are equivalent :
\begin{enumerate}
\item There is some $a>0$ such that $\mu$ satisfies $\ST_\a(a)$.
\item There is some $b>0$ such that
\[K^+:=\sup_{x\geq m} \int e^{\a(bz)}\,d\mu_x^+(z)<+\infty\quad\text{and}\quad K^-:=\sup_{x\leq m}\int e^{\a(bz)}\,d\mu_x^-(z)<+\infty.\]
\end{enumerate}
\end{thm}
\proof~\\
$[(1)\Rightarrow(2)]$ According to Proposition \ref{integrability2}, if $\mu$ satisfies $\ST_\a(a)$, then
\begin{align*}
&\int e^{\a(bz)}\,d\mu_x^+(z)\leq \frac{1}{\mu (-\infty,x]}+1\leq 3,\quad \forall x\geq m\\
&\int e^{\a(bz)}\,d\mu_x^-(z)\leq \frac{1}{\mu [x,+\infty)}+1\leq 3,\quad \forall x\leq m.
\end{align*}
Thus (2) holds true for $b=a$.

$[(2)\Rightarrow(1)]$
According to Corollary \ref{unif.stoc.dom}, if there is some $a>0$ such that
\begin{align}
&\sup_{x\geq m} \mu_x^+[h,+\infty)\leq e^{-\a_1^{-1}\circ\a(ah)},\quad \forall h\geq 0\label{char.lm.1}\\
&\sup_{x\leq m} \mu_x^-[h,+\infty)\leq e^{-\a_1^{-1}\circ\a(ah)},\quad \forall h\geq 0\label{char.lm.2},
\end{align}
then $\mu$ satisfies the strong TCI with the cost function $\frac 1 \kappa \a(a|x-y|/2)\geq \a(a|x-y|/(2\kappa))$. Hence it is enough to prove (\ref{char.lm.1}) and (\ref{char.lm.2}).

Let us prove (\ref{char.lm.1}) (the proof of (\ref{char.lm.2}) will be the same). To prove (\ref{char.lm.1}), it is enough to find $a>0$ such that
\begin{equation}\label{char.lm.3}
\la\begin{array}{ll}\mu_x^+[h,+\infty)\leq e^{-ah},&\quad \forall h\leq \frac 1 a\\
\mu_x^+[h,+\infty)\leq e^{-\a(ah)},&\quad \forall h\geq \frac 1 a
\end{array}\right.,
\end{equation}
holds for all $x\geq m$.
Since $\mu\in \lm$, there is $a_0>0$ such that
\[\sup_{x\geq m}\mu_x^+[h,+\infty) \leq e^{-a_0 h},\quad \forall h\geq 0,\]
and by hypothesis there is some $b_0>0$ such that $K^+:=\sup_{x\geq m} \int e^{\a(b_0z)}\,d\mu_x^+(z)<+\infty.$
To conclude it suffices to apply Lemma \ref{tech.lem} with $\nu=\mu_x^+$ and $K=K^+$. It provides us a constant $a>0$ depending only on $a_0,b_0$ and $K^+$ such that (\ref{char.lm.3}) holds true for all $x\geq m$.
\endproof

\subsection{Tractable sufficient condition for good probability measures}
\begin{thm}\label{suff.cond}
Let $d\mu=e^{-V}\,dx$ with $V\in \V$ and $\a\in \A\cap \V$. If $\mu\in \lm$ and if
\begin{equation}\label{a'/V'}
\exists \lambda>0 \text{ such that }\limsup_{u\to\pm \infty}\frac{\a'(\lambda u)}{V'(u+m)}<+\infty,
\end{equation}
where $m$ is the median of $\mu$,
then $\mu$ satisfies the inequality $\ST_\a(a)$ for some $a>0$.
\end{thm}
To prove this theorem, we will use the following Lemma.
\begin{lem}\label{int.equiv}
Let $\Phi \in \V$, then \[\int_x^{+\infty}e^{-\Phi(t)}\,dt\sim \frac{e^{-\Phi(x)}}{\Phi'(x)}\quad\text{and}\quad\int_0^{x}e^{-\Phi(t)}\,dt\sim \frac{e^{\Phi(x)}}{\Phi'(x)},\quad\text{as }x \text{ goes to }+\infty.\]
\end{lem}
\proof See Corollary 6.4.2 of \cite{Log-Sob}. \endproof
\proof[Proof of Theorem \ref{suff.cond}]
Let $\tilde{\mu}=\mathcal{L}(X-m)$, where $X$ is a random measure with law $\mu$. The density of $\tilde{\mu}$ with respect to Lebesgues measure is $e^{-\widetilde{V}}$, with $\widetilde{V}(x)=V(x+m),\forall x\in \R$. As $x\mapsto x+m$ is $1$-Lipschitz, it follows from Corollary \ref{contr.princ} that $\mu$ satisfies $\ST_\a(a)$ if and only if $\tilde{\mu}$ satisfies $\ST_\a(a)$.
Observe that $\tilde{\mu}\in \lm$. According to Theorem \ref{char.lm}, to prove that $\tilde{\mu}$ satisfies $\ST_\a(a)$ for some $a>0$, it suffices to prove that there is $b>0$ such that
\[K^+(b)=\sup_{x\geq 0}\int e^{\a(bu)}\,d\tilde{\mu}^+_x(u)<+\infty\qquad\text{and}\qquad K^-(b)=\sup_{x\leq 0}\int e^{\a(bu)}\,d\tilde{\mu}^-_x(u)<+\infty,\]
where $\tilde{\mu}^+_x=\mathcal{L}(\widetilde{X}-x|\widetilde{X}\geq x)$ and $\tilde{\mu}^-_x=\mathcal{L}(x-\widetilde{X}|\widetilde{X}\leq x)$
with $\widetilde{X}$ of law $\tilde{\mu}$.

The proof of $K^-(b)<+\infty$ being similar, we will only prove that $K^+(b)<+\infty$ for some $b>0$. One can suppose without restriction that $\lambda=1$ in (\ref{a'/V'}).
Define
\begin{align*}
K(b,x)&=\int e^{ \a(bt)}\,d\tilde{\mu}_x^+(t)\\
& = \frac{\int_x^{+\infty}e^{\a(b(u-x))}e^{-\widetilde{V}(u)}\,du}{\int_x^{+\infty} e^{-\widetilde{V}(u)}\,du},\qquad \forall x\geq 0,\forall b\geq 0.
\end{align*}
Let us show that there is $k\in \N^*$ such that $K(1/k,x)<+\infty$ for all $x\geq0$.
Since $\a$ is super-additive and non decreasing, one gets
\[K(1/k,x)\leq \frac{\int_0^{+\infty}e^{\frac 1 k\a(u)} e^{-\widetilde{V}(u)}\,du}{\int_x^{+\infty}e^{-\widetilde{V}(u)}\,du}.\]
Since $\displaystyle{\limsup_{u\rightarrow+\infty}\frac{\a'(u)}{\widetilde{V}'(u)}<+\infty}$, there are $M>0$ and $u_0>0$ such that
$\widetilde{V}'(u)\geq M \a'(u),$ for all $u\geq u_0$.
Integrating yields
\[\widetilde {V}(u)\geq M \a(u)+C,\quad \forall u\geq u_0,\]
where $C$ is a constant.
Let $k_0$ be a positive integer such that $k_0\geq \frac 2 M$. Then one has
\[
e^{\frac{1}{k_0}\a(u)-\widetilde{V}(u)} \leq e^{-\frac{M}{2}\a(u)-C}
 \leq e^{-\frac{M\a(1)}{2}(u-1)-C},
\quad \forall u\geq u_0\]
where the last inequality follows from the inequality $\a(u)\geq \a(1)(u-1),\forall u\geq 0$ which is easy to establish. From this follows easily that $K(1/k_0,x)<+\infty$ for all $x\geq0$.

Now, let us show that $\sup_{x\geq 0} K(1/k_0,x)<+\infty.$ Since the map $x\mapsto K(1/k_0,x)$ is continuous, it suffices to check that $\displaystyle{\limsup_{x\to+\infty}K(1/k_0,x)<+\infty}.$
Using the super-additivity of $\a$, one gets
\[\a(u-x)\leq \a(u)-\a(x),\quad\forall u\geq x\geq 0.\]
So
\[K(1/k_0,x)\leq e^{-1/k_0\a(x)}\frac{\int_x^{+\infty}e^{1/k_0\a(u)-\widetilde{V}(u)}\,du}{\int_x^{+\infty}e^{-\widetilde{V}(u)}\,du}.\]
Applying Lemma \ref{int.equiv}, with $\Phi=\widetilde{V}-1/k_0\a$, and then with $\Phi=\widetilde{V}$, one gets
\[e^{-1/k_0\a(x)}\frac{\int_x^{+\infty}e^{1/k_0\a(u)-\widetilde{V}(u)}\,du}{\int_x^{+\infty}e^{-\widetilde{V}(u)}\,du}\thicksim e^{-1/k_0\a(x)}\frac{e^{1/k_0\a(x)-\widetilde{V}(x)}}{\widetilde{V}'(x)-1/k_0\a'(x)}e^{\widetilde{V}(x)}\widetilde{V}'(x)=\frac{1}{1-1/k_0\frac{\a'(x)}{\widetilde{V}'(x)}}.\]
Since $\displaystyle{\limsup_{x\to +\infty} \frac{1}{1-1/k_0\frac{\a'(x)}{\widetilde{V}'(x)}}<+\infty}$, one deduces that $\displaystyle{\limsup_{x\to+\infty}K(1/k_0,x)<+\infty},$ which ends the proof.
\endproof
\bibliographystyle{plain}
\bibliography{bib}

\end{document}